\documentclass[10pt]{article}
\usepackage{cite}
\usepackage{mathrsfs}
\usepackage{amsfonts}
\usepackage{amsmath}
\usepackage{amsfonts,amssymb}
\usepackage{dsfont}
\usepackage{curves}
\usepackage{mathrsfs}
\usepackage{pifont}
\usepackage{amssymb}
\allowdisplaybreaks

\numberwithin{equation}{section}

\date{}

\def\BigRoman{\uppercase\expandafter{\romannumeral\number\count 255 }}
\def\Romannumeral{\afterassignment\BigRoman\count255=}

\begin{document}
\title{Toughness and $A_{\alpha}$-spectral radius in graphs
}
\author{\small  Sizhong Zhou$^{1}$\footnote{Corresponding
author. E-mail address: zsz\_cumt@163.com (S. Zhou)}, Yuli Zhang$^{2}$\footnote{Corresponding
author. E-mail address: zhangyuli\_djtu@126.com (Y. Zhang)}, Tao Zhang$^{3}$, Hongxia Liu$^{4}$\\
\small $1$. School of Science, Jiangsu University of Science and Technology,\\
\small Zhenjiang, Jiangsu 212100, China\\
\small $2$. School of Science, Dalian Jiaotong University,\\
\small Dalian, Liaoning 116028, China\\
\small $3$. School of Economics and management, Jiangsu University of Science and Technology,\\
\small Zhenjiang, Jiangsu 212100, China\\
\small $4$. School of Mathematics and Information Science, Yantai University,\\
\small Yantai, Shandong 264005, China\\
}

\maketitle
\begin{abstract}
\noindent Let $\alpha\in[0,1)$, and let $G$ be a connected graph of order $n$ with $n\geq f(\alpha)$, where $f(\alpha)=6$ for $\alpha\in[0,\frac{2}{3}]$
and $f(\alpha)=\frac{4}{1-\alpha}$ for $\alpha\in(\frac{2}{3},1)$. A graph $G$ is said to be $t$-tough if $|S|\geq tc(G-S)$ for each subset $S$ of
$V(G)$ with $c(G-S)\geq2$, where $c(G-S)$ is the number of connected components in $G-S$. The $A_{\alpha}$-spectral radius of $G$ is denoted by
$\rho_{\alpha}(G)$. In this paper, it is verified that $G$ is a 1-tough graph unless $G=K_1\vee(K_{n-2}\cup K_1)$ if $\rho_{\alpha}(G)\geq\rho_{\alpha}(K_1\vee(K_{n-2}\cup K_1))$, where $\rho_{\alpha}(K_1\vee(K_{n-2}\cup K_1))$ equals the largest root of $x^{3}-((\alpha+1)n+\alpha-3)x^{2}+(\alpha n^{2}+(\alpha^{2}-\alpha-1)n-2\alpha+1)x-\alpha^{2}n^{2}+(3\alpha^{2}-\alpha+1)n-4\alpha^{2}+5\alpha-3=0$.
Further, we present an $A_{\alpha}$-spectral radius condition for a graph to be a $t$-tough graph.
\\
\begin{flushleft}
{\em Keywords:} graph; toughness; $A_{\alpha}$-spectral radius.

(2020) Mathematics Subject Classification: 05C50
\end{flushleft}
\end{abstract}

\section{Introduction}

Let $G$ be an undirected simple graph with vertex set $V(G)=\{v_1,v_2,\ldots,v_n\}$ and edge set $E(G)$. The order of $G$ is the number
$n=|V(G)|$ of its vertices and its size is the number $m=e(G)$ of its edges. A graph $G$ is called trivial if $|V(G)|=1$. For arbitrary
$v\in V(G)$, the degree of $v$ in $G$ is defined as
the number of edges which are adjacent to $v$ and denoted by $d_G(v)$. For a given subset $S$ of $V(G)$, we use $G[S]$ to denote the subgraph
of $G$ induced by $S$, and write $G-S$ for $G[V(G)\setminus S]$. The complement of a graph $G$ is a graph $\overline{G}$ with the same vertex
set as $G$, in which any two distinct vertices are adjacent if and only if they are nonadjacent in $G$. Let $G_1$ and $G_2$ be two
vertex-disjoint graphs. We use $G_1\cup G_2$ to denote the disjoint union of $G_1$ and $G_2$. The join $G_1\vee G_2$ is the graph formed
from $G_1\cup G_2$ by adding all possible edges between $V(G_1)$ and $V(G_2)$. We denote by $K_n$ a complete graph of order $n$.

Let $t$ be a nonnegative real number. A graph $G$ is said to be $t$-tough if $|S|\geq tc(G-S)$ for each subset $S$ of $V(G)$ with
$c(G-S)\geq2$, where $c(G-S)$ is the number of connected components in $G-S$. The toughness $t(G)$ of $G$ is the largest real number $t$ for
which $G$ is $t$-tough, or is $\infty$ if $G$ is complete. This concept was first introduced by Chv\'atal \cite{C} in 1973. Some results on
toughness of graphs can be found in \cite{CGL,CLW,EJKS,LFS,WZi,Zhs,Zt,ZXS,ZWB,ZWX}.

Given a graph $G$ of order $n$, the adjacency matrix $A(G)$ of $G$ is the $n\times n$ matrix in which entry $a_{ij}$ is 1 or 0 according to
whether $v_i$ and $v_j$ are adjacent or not, and $a_{ij}=0$ if $i=j$. The eigenvalues of $G$ are the eigenvalues of its adjacency matrix $A(G)$. Let
$\lambda_1(G)\geq\lambda_2(G)\geq\dots\geq\lambda_n(G)$ be its eigenvalues in nonincreasing order. Note that the adjacency spectral radius
of $G$ is equal to $\lambda_1(G)$, written as $\rho(G)$.

Let $D(G)$ denote the diagonal matrix of vertex degree of $G$. Then $L(G)=D(G)-A(G)$ and
$Q(G)=D(G)+A(G)$ are called the Laplacian matrix and signless Laplacian matrix of $G$, respectively. For any $\alpha\in[0,1]$, Nikiforov
\cite{N} introduced the $A_{\alpha}$-matrix of $G$ as
$$
A_{\alpha}(G)=\alpha D(G)+(1-\alpha)A(G).
$$
It is easy to see that $A_{\alpha}(G)=A(G)$ if $\alpha=0$, and $A_{\alpha}(G)=\frac{1}{2}Q(G)$ if $\alpha=\frac{1}{2}$. Note that
$A_{\alpha}(G)$ is a real symmetric nonnegative matrix. Hence, the eigenvalues of $A_{\alpha}(G)$ are real, which can be indexed in
nonincreasing order as $\lambda_1(A_{\alpha}(G))\geq\lambda_2(A_{\alpha}(G))\geq\cdots\geq\lambda_n(A_{\alpha}(G))$. The largest eigenvalue
$\lambda_1(A_{\alpha}(G))$ is called the $A_{\alpha}$-spectral radius of $G$, and denoted by $\rho_{\alpha}(G)$. Namely,
$\rho_{\alpha}(G)=\lambda_1(A_{\alpha}(G))$. For some interesting spectral properties of $A_{\alpha}(G)$, we refer the reader to
\cite{BFO,HLZ,LHX,LLX,LXS,N,NR,ZHW}.

Many authors \cite{FLL,O,Wc,Ws,Zs,ZW,ZW1,ZZ,ZZL} obtained some results on spectral radius in graphs. Brouwer \cite{B} investigated the
relationship between toughness and eigenvalues and claimed that for any connected $d$-regular graph $G$,
$t(G)>\frac{d}{\lambda}-2$, and he further conjectured that the lower bound $\frac{d}{\lambda}-2$ can be improved to $\frac{d}{\lambda}-1$,
where $\lambda=\max_{2\leq i\leq n}|\lambda_i|$. Later, Gu \cite{Gt} strengthened Brouwer's result and claimed the lower bound
$\frac{d}{\lambda}-2$ can be improved to $\frac{d}{\lambda}-\sqrt{2}$. Very recently, Gu \cite{Ga} confirmed Brouwer's conjecture. Fan, Lin
and Lu \cite{FLLt} presented two adjacency spectral radius conditions for the existence of a $t$-tough graph. It is natural and interesting
to find a sufficient condition for a graph to be $t$-tough in view of its $A_{\alpha}$-spectral radius. In this paper, we first show a
$A_{\alpha}$-spectral radius condition to guarantee that a graph is 1-tough.

\medskip

\noindent{\textbf{Theorem 1.1.}} Let $\alpha\in[0,1)$, and let $G$ be a connected graph of order $n$ with $n\geq f(\alpha)$, where
\[
f(\alpha)=\left\{
\begin{array}{ll}
6,&if \ \alpha\in[0,\frac{2}{3}];\\
\frac{4}{1-\alpha},&if \ \alpha\in(\frac{2}{3},1).\\
\end{array}
\right.
\]
If $\rho_{\alpha}(G)\geq\rho_{\alpha}(K_1\vee(K_{n-2}\cup K_1))$, then $G$ is a 1-tough graph unless $G=K_1\vee(K_{n-2}\cup K_1)$, where
$\rho_{\alpha}(K_1\vee(K_{n-2}\cup K_1))$ equals the largest root of $x^{3}-((\alpha+1)n+\alpha-3)x^{2}+(\alpha n^{2}+(\alpha^{2}-\alpha-1)n-2\alpha+1)x-\alpha^{2}n^{2}+(3\alpha^{2}-\alpha+1)n-4\alpha^{2}+5\alpha-3=0$.

\medskip

Further, it is natural and interesting to find an $A_{\alpha}$-spectral radius condition to ensure that a graph is a $t$-tough graph. Next, we
present an $A_{\alpha}$-spectral radius condition for a graph to be $t$-tough.

\medskip

\noindent{\textbf{Theorem 1.2.}} Let $\alpha\in[\frac{1}{2},\frac{3}{4})$, and let $t$ be a positive integer. If $G$ is a connected graph of order
$n\geq\max\{5t^{2}+10t+1,\frac{12t(1-\alpha)-2\alpha+1}{3-4\alpha}\}$ with $\rho_{\alpha}(G)\geq\rho_{\alpha}(K_{2t-1}\vee(K_{n-2t}\cup K_1))$, then
$G$ is a $t$-tough graph unless $G=K_{2t-1}\vee(K_{n-2t}\cup K_1)$.

\section{Preliminaries}

In this section, we put forward some necessary lemmas, which play a key role in proving our main results.

\medskip

\noindent{\textbf{Lemma 2.1}} (Nikiforov \cite{N}). Let $K_n$ be a complete graph of order $n$. Then
$$
\rho_{\alpha}(K_n)=n-1.
$$

\medskip

\noindent{\textbf{Lemma 2.2}} (Nikiforov \cite{N}). If $G$ is a connected graph, and $H$ is a proper subgraph of $G$, then
$$
\rho_{\alpha}(G)>\rho_{\alpha}(H).
$$

\medskip

\noindent{\textbf{Lemma 2.3}} (Zhao, Huang and Wang \cite{ZHW}). Let $\alpha\in[0,1)$, and let $n_1\geq n_2\geq\cdots\geq n_t$ be positive
integers with $n=\sum_{i=1}^{t}n_i+s$ and $n_1\leq n-s-t+1$. Then
$$
\rho_{\alpha}(K_s\vee(K_{n_1}\cup K_{n_2}\cup\cdots\cup K_{n_t}))\leq\rho_{\alpha}(K_s\vee(K_{n-s-t+1}\cup(t-1)K_1)),
$$
where the equality holds if and only if $(n_1,n_2,\ldots,n_t)=(n-s-t+1,1,\ldots,1)$.

\medskip

Let $M$ be a real matrix whose rows and columns are indexed by $V=\{1,2,\ldots,n\}$. Assume that $M$, with respect to the partition
$\pi:V=V_1\cup V_2\cup\cdots\cup V_t$, can be written as
\begin{align*}
M=\left(
  \begin{array}{cccc}
    M_{11} & M_{12} & \cdots & M_{1t}\\
    M_{21} & M_{22} & \cdots & M_{2t}\\
    \vdots & \vdots & \ddots & \vdots\\
    M_{t1} & M_{t2} & \cdots & M_{tt}\\
  \end{array}
\right),
\end{align*}
where $M_{ij}$ denotes the submatrix (block) of $M$ formed by rows in $V_i$ and columns in $V_j$. Let $b_{ij}$ denote the average row sum of
$M_{ij}$, namely, $b_{ij}$ is the sum of all entries in $M_{ij}$ divided by the number of rows. Then matrix $M_{\pi}=(b_{ij})$ is called the
\emph{quotient matrix} of $M$. If the row sum of every block $M_{ij}$ is a constant, then the partition is \emph{equitable}.

\medskip

\noindent{\textbf{Lemma 2.4}} (Haemers \cite{H}, You, Yang, So and Xi \cite{YYSX}). Let $M$ be a real matrix with an equitable partition $\pi$,
and let $M_{\pi}$ be the corresponding quotient matrix. Then each eigenvalue of $M_{\pi}$ is an eigenvalue of $M$. Furthermore, if $M$ is
nonnegative, then the spectral radius of $M_{\pi}$ equals the spectral radius of $M$.

\medskip

\noindent{\textbf{Lemma 2.5}} (Alhevaz, Baghipur, Ganie and Das \cite{ABGD}). Let $G$ be a graph of order $n$ with $m$ edges, with no isolated
vertices and let $\alpha\in[\frac{1}{2},1]$. Then
$$
\rho_{\alpha}(G)\leq\frac{2m(1-\alpha)}{n-1}+\alpha n-1.
$$
If $\alpha\in(\frac{1}{2},1)$ and $G$ is connected, equality holds if and only if $G=K_n$.

\section{The proof of Theorem 1.1}

In this section, we verify Theorem 1.1, which establishes an $A_{\alpha}$-spectral radius condition for a graph to be 1-tough.

\medskip

\noindent{\it Proof of Theorem 1.1.} Suppose to the contrary that $G$ is not a 1-tough graph, there exists some nonempty subset $S\subseteq V(G)$
such that $c(G-S)\geq|S|+1$. Let $|S|=s$. Then $G$ is a spanning subgraph of $G_s^{1}=K_s\vee(K_{n_1}\cup K_{n_2}\cup\cdots\cup K_{n_{s+1}})$ for
some positive integers $n_1\geq n_2\geq\cdots\geq n_{s+1}$ with $\sum\limits_{i=1}^{s+1}n_i=n-s$. According to Lemma 2.2, we have
\begin{align}\label{eq:3.1}
\rho_{\alpha}(G)\leq\rho_{\alpha}(G_s^{1})
\end{align}
with equality if and only if $G=G_s^{1}$.

Let $G_{s}^{2}=K_s\vee(K_{n-2s}\cup sK_1)$. Using Lemma 2.3, we get
\begin{align}\label{eq:3.2}
\rho_{\alpha}(G_s^{1})\leq\rho_{\alpha}(G_s^{2}),
\end{align}
where the equality holds if and only if $(n_1,n_2,\ldots,n_{s+1})=(n-2s,1,\ldots,1)$. Let $\varphi(x)=x^{3}-((\alpha+1)n+\alpha-3)x^{2}+(\alpha n^{2}+(\alpha^{2}-\alpha-1)n-2\alpha+1)x-\alpha^{2}n^{2}+(3\alpha^{2}-\alpha+1)n-4\alpha^{2}+5\alpha-3$. Then $\rho_{\alpha}(K_1\vee(K_{n-2}\cup K_1))$
equals the largest root of $\varphi(x)=0$. We are to verify the following claim.

\medskip

\noindent{\bf Claim 1.} If $\alpha\in[0,1)$ and $n\geq f(\alpha)$, then we obtain
$$
\rho_{\alpha}(G_s^{2})\leq\rho_{\alpha}(K_1\vee(K_{n-2}\cup K_1))
$$
with equality if and only if $G_s^{2}=K_1\vee(K_{n-2}\cup K_1)$.

\noindent{\it Proof.} Recall that $G_{s}^{2}=K_s\vee(K_{n-2s}\cup sK_1)$. Obviously, $1\leq s\leq\frac{n-1}{2}$.

Consider the partition $V(G_{s}^{2})=V(sK_1)\cup V(K_{n-2s})\cup V(K_s)$. The corresponding quotient matrix of $A_{\alpha}(G_{s}^{2})$ is
equal to
\begin{align*}
B_1=\left(
  \begin{array}{ccc}
    \alpha s & 0 & (1-\alpha)s\\
    0 & n+(\alpha-2)s-1 & (1-\alpha)s\\
    (1-\alpha)s & (1-\alpha)(n-2s) & \alpha n-\alpha s+s-1\\
  \end{array}
\right).
\end{align*}

Then the characteristic polynomial of $B_1$ is
\begin{align}\label{eq:3.3}
\varphi_{B_1}(x)=&x^{3}-((\alpha+1)n+(\alpha-1)s-2)x^{2}\nonumber\\
&+(\alpha n^{2}+(\alpha^{2}s-\alpha-1)n-s^{2}-(2\alpha-1)s+1)x\nonumber\\
&-\alpha^{2}sn^{2}+(2\alpha^{2}-2\alpha+1)s^{2}n+(\alpha^{2}+\alpha)sn\nonumber\\
&-(3\alpha^{2}-5\alpha+2)s^{3}-(\alpha^{2}-\alpha+1)s^{2}-\alpha s.
\end{align}
Since the partition $V(G_{s}^{2})=V(sK_1)\cup V(K_{n-2s})\cup V(K_s)$ is equitable, it follows from Lemma 2.4 that $\rho_{\alpha}(G_{s}^{2})$
is the largest root of $\varphi_{B_1}(x)=0$. Namely, $\varphi_{B_1}(\rho_{\alpha}(G_{s}^{2}))=0$. Let
$\eta_1=\rho_{\alpha}(G_{s}^{2})\geq\eta_2\geq\eta_3$ be the three roots of $\varphi_{B_1}(x)=0$ and $Q=\mbox{diag}(s,n-2s,s)$. One checks that
$$
Q^{\frac{1}{2}}B_1Q^{-\frac{1}{2}}=\left(
  \begin{array}{ccc}
    \alpha s & 0 & (1-\alpha)s\\
    0 & n+(\alpha-2)s-1 & (1-\alpha)s^{\frac{1}{2}}(n-2s)^{\frac{1}{2}}\\
    (1-\alpha)s & (1-\alpha)s^{\frac{1}{2}}(n-2s)^{\frac{1}{2}} & \alpha n-\alpha s+s-1\\
  \end{array}
\right)
$$
is symmetric, and also contains
\begin{align*}
\left(
  \begin{array}{ccc}
    \alpha s & 0\\
    0 & n+(\alpha-2)s-1\\
  \end{array}
\right)
\end{align*}
as its submatrix. Since $Q^{\frac{1}{2}}B_1Q^{-\frac{1}{2}}$ and $B_1$ have the same eigenvalues, the Cauchy interlacing theorem (cf. \cite{H})
implies that
\begin{align}\label{eq:3.4}
\eta_2\leq n+(\alpha-2)s-1<n-2.
\end{align}

If $s=1$, then $G_s^{2}=K_1\vee(K_{n-2}\cup K_1)$ and $\varphi_{B_1}(x)=\varphi(x)$. And so $\rho_{\alpha}(G_s^{2})=\rho_{\alpha}(K_1\vee(K_{n-2}\cup K_1))$.
In what follows, we always assume $s\geq2$.

Note that $K_{n-1}$ is a proper subgraph of $K_1\vee(K_{n-2}\cup K_1)$. In terms of \eqref{eq:3.4}, Lemmas 2.1 and 2.2, we get
\begin{align}\label{eq:3.5}
\rho_{\alpha}(K_1\vee(K_{n-2}\cup K_1))>\rho_{\alpha}(K_{n-1})=n-2>\eta_2.
\end{align}

Write $\theta(n)=\rho_{\alpha}(K_1\vee(K_{n-2}\cup K_1))$. Note that $\varphi(\theta(n))=0$. A simple calculation yields that
\begin{align}\label{eq:3.6}
\varphi_{B_1}(\theta(n))=&\varphi_{B_1}(\theta(n))-\varphi(\theta(n))\nonumber\\
=&(s-1)((1-\alpha)(\theta(n))^{2}+(\alpha^{2}n-2\alpha-s)\theta(n)-\alpha^{2}n^{2}\nonumber\\
&+(2\alpha^{2}-2\alpha+1)sn+(3\alpha^{2}-\alpha+1)n-(3\alpha^{2}-5\alpha+2)s^{2}\nonumber\\
&-(4\alpha^{2}-6\alpha+3)s-4\alpha^{2}+5\alpha-3).
\end{align}

The following proof will be divided into two cases.

\noindent{\bf Case 1.} $0\leq\alpha\leq\frac{2}{3}$.

Let $h(x)=(1-\alpha)x^{2}+(\alpha^{2}n-2\alpha-s)x-\alpha^{2}n^{2}+(2\alpha^{2}-2\alpha+1)sn+(3\alpha^{2}-\alpha+1)n-(3\alpha^{2}-5\alpha+2)s^{2}
-(4\alpha^{2}-6\alpha+3)s-4\alpha^{2}+5\alpha-3$ be a real function in $x$. Then the symmetry axis of $h(x)$ is
$x=\frac{-\alpha^{2}n+2\alpha+s}{2(1-\alpha)}$. Obviously, $h(x)$ is increasing when $x\geq\frac{-\alpha^{2}n+2\alpha+s}{2(1-\alpha)}$.
By plugging the value $\theta(n)$ into $x$ of $h(x)$, we obtain
\begin{align*}
h(\theta(n))=&(1-\alpha)(\theta(n))^{2}+(\alpha^{2}n-2\alpha-s)\theta(n)-\alpha^{2}n^{2}\\
&+(2\alpha^{2}-2\alpha+1)sn+(3\alpha^{2}-\alpha+1)n-(3\alpha^{2}-5\alpha+2)s^{2}\\
&-(4\alpha^{2}-6\alpha+3)s-4\alpha^{2}+5\alpha-3
\end{align*}
and it follows from \eqref{eq:3.6} that
\begin{align}\label{eq:3.7}
\varphi_{B_1}(\theta(n))=(s-1)h(\theta(n)).
\end{align}
Note that $n-2s\geq1$. We have $n\geq2s+1$. Together with $s\geq2$ and \eqref{eq:3.5}, we conclude
$$
\frac{-\alpha^{2}n+2\alpha+s}{2(1-\alpha)}<n-2<\rho_{\alpha}(K_1\vee(K_{n-2}\cup K_1))=\theta(n),
$$
and so
\begin{align}\label{eq:3.8}
h(\theta(n))>&h(n-2)\nonumber\\
=&(1-\alpha)(n-2)^{2}+(\alpha^{2}n-2\alpha-s)(n-2)-\alpha^{2}n^{2}\nonumber\\
&+(2\alpha^{2}-2\alpha+1)sn+(3\alpha^{2}-\alpha+1)n-(3\alpha^{2}-5\alpha+2)s^{2}\nonumber\\
&-(4\alpha^{2}-6\alpha+3)s-4\alpha^{2}+5\alpha-3\nonumber\\
=&(1-\alpha)n^{2}+((2\alpha^{2}-2\alpha)s+\alpha^{2}+\alpha-3)n-(3\alpha^{2}-5\alpha+2)s^{2}\nonumber\\
&-(4\alpha^{2}-6\alpha+1)s-4\alpha^{2}+5\alpha+1.
\end{align}
Let $p(x,s)=(1-\alpha)x^{2}+((2\alpha^{2}-2\alpha)s+\alpha^{2}+\alpha-3)x-(3\alpha^{2}-5\alpha+2)s^{2}-(4\alpha^{2}-6\alpha+1)s-4\alpha^{2}+5\alpha+1$
be a real function in $x$. According to \eqref{eq:3.8}, we get
\begin{align}\label{eq:3.9}
h(\theta(n))>p(n,s).
\end{align}
Clearly, $p(x,s)$ is increasing when $x\geq-\frac{(2\alpha^{2}-2\alpha)s+\alpha^{2}+\alpha-3}{2(1-\alpha)}$. If $s\geq3$, then we can verify that
$$
-\frac{(2\alpha^{2}-2\alpha)s+\alpha^{2}+\alpha-3}{2(1-\alpha)}<2s+1\leq n,
$$
and so
\begin{align*}
p(n,s)\geq&p(2s+1,s)\\
=&(1-\alpha)(2s+1)^{2}+((2\alpha^{2}-2\alpha)s+\alpha^{2}+\alpha-3)(2s+1)\\
&-(3\alpha^{2}-5\alpha+2)s^{2}-(4\alpha^{2}-6\alpha+1)s-4\alpha^{2}+5\alpha+1\\
=&(s^{2}-3)\alpha^{2}-(3s^{2}-2s-5)\alpha+2s^{2}-3s-1\\
\geq&\frac{4}{9}(s^{2}-3)-\frac{2}{3}(3s^{2}-2s-5)+2s^{2}-3s-1\\
=&\frac{4s^{2}-15s+9}{9}\\
\geq&0,
\end{align*}
where the last two inequalities hold from $\frac{3s^{2}-2s-5}{2(s^{2}-3)}>\frac{2}{3}\geq\alpha$ and $s\geq3$, respectively.

If $s=2$, then $-\frac{(2\alpha^{2}-2\alpha)s+\alpha^{2}+\alpha-3}{2(1-\alpha)}=-\frac{5\alpha^{2}-3\alpha-3}{2(1-\alpha)}<6=f(\alpha)\leq n$
due to $0\leq\alpha\leq\frac{2}{3}$. Hence, for $0\leq\alpha\leq\frac{2}{3}$ and $n\geq f(\alpha)=6$, we get
\begin{align*}
p(n,2)\geq p(6,2)=6\alpha^{2}-17\alpha+9>0.
\end{align*}

Thus, we obtain $p(n,s)>0$ for $2\leq s\leq\frac{n-1}{2}$. Together with \eqref{eq:3.7}, \eqref{eq:3.9} and $s\geq2$, we have
$$
\varphi_{B_1}(\theta(n))=(s-1)h(\theta(n))>(s-1)p(n,s)\geq0.
$$
Recall that $\rho_{\alpha}(G_{s}^{2})$ is the largest root of $\varphi_{B_1}(x)=0$. From \eqref{eq:3.6}, the symmetry axis of $\varphi_{B_1}(\theta(n))$
is $\theta(n)=\frac{-\alpha^{2}n+2\alpha+s}{2(1-\alpha)}$. Thus, we see that $\varphi_{B_1}(\theta(n))$ is increasing in the interval
$[\frac{-\alpha^{2}n+2\alpha+s}{2(1-\alpha)},+\infty)$. Note that $\frac{-\alpha^{2}n+2\alpha+s}{2(1-\alpha)}<n-2$. Together with \eqref{eq:3.5}, we
have $\frac{-\alpha^{2}n+2\alpha+s}{2(1-\alpha)}<n-2<\theta(n)$. Combining these with $\varphi_{B_1}(\theta(n))>0$, we deduce
$$
\rho_{\alpha}(G_{s}^{2})<\theta(n)=\rho_{\alpha}(K_1\vee(K_{n-2}\cup K_1))
$$
for $2\leq s\leq\frac{n-1}{2}$.

\noindent{\bf Case 2.} $\frac{2}{3}<\alpha<1$.

In terms of \eqref{eq:3.3}, we obtain
\begin{align*}
\varphi_{B_1}(n-2)=&(n-2)^{3}-((\alpha+1)n+(\alpha-1)s-2)(n-2)^{2}\\
&+(\alpha n^{2}+(\alpha^{2}s-\alpha-1)n-s^{2}-(2\alpha-1)s+1)(n-2)\\
&-\alpha^{2}sn^{2}+(2\alpha^{2}-2\alpha+1)s^{2}n+(\alpha^{2}+\alpha)sn\\
&-(3\alpha^{2}-5\alpha+2)s^{3}-(\alpha^{2}-\alpha+1)s^{2}-\alpha s\\
=&(3\alpha-2)(1-\alpha)s^{3}+((2\alpha^{2}-2\alpha)n-\alpha^{2}+\alpha+1)s^{2}\\
&+((1-\alpha)n^{2}-(\alpha^{2}-3\alpha+3)n-\alpha+2)s\\
&+(\alpha-1)n^{2}-(2\alpha-3)n-2\\
&:=\Phi(s,n).
\end{align*}
Thus, we obtain
\begin{align*}
\frac{\partial \Phi(s,n)}{\partial s}=&3(3\alpha-2)(1-\alpha)s^{2}+2((2\alpha^{2}-2\alpha)n-\alpha^{2}+\alpha+1)s\\
&+(1-\alpha)n^{2}-(\alpha^{2}-3\alpha+3)n-\alpha+2.
\end{align*}
Note that $\frac{2}{3}<\alpha<1$ and $n\geq f(\alpha)=\frac{4}{1-\alpha}$. By a simple calculation, we have
\begin{align*}
\frac{\partial \Phi(s,n)}{\partial s}\Big|_{s=2}=&(1-\alpha)n^{2}+(7\alpha^{2}-5\alpha-3)n-40\alpha^{2}+63\alpha-18\\
\geq&(1-\alpha)\left(\frac{4}{1-\alpha}\right)^{2}+(7\alpha^{2}-5\alpha-3)\left(\frac{4}{1-\alpha}\right)-40\alpha^{2}+63\alpha-18\\
=&\frac{1}{1-\alpha}(40\alpha^{3}-75\alpha^{2}+61\alpha-14)\\
>&0,
\end{align*}
and
\begin{align*}
\frac{\partial \Phi(s,n)}{\partial s}\Big|_{s=\frac{n-1}{2}}=&\frac{1}{4}((1-\alpha)(\alpha-2)n^{2}+(2\alpha^{2}-6\alpha-4)n-5\alpha^{2}+7\alpha+6)\\
\leq&\frac{1}{4}\left((1-\alpha)(\alpha-2)\left(\frac{4}{1-\alpha}\right)^{2}+(2\alpha^{2}-6\alpha-4)\left(\frac{4}{1-\alpha}\right)-5\alpha^{2}+7\alpha+6\right)\\
=&\frac{1}{4(1-\alpha)}(5\alpha^{3}-4\alpha^{2}-7\alpha-42)\\
<&0.
\end{align*}
This implies that $\varphi_{B_1}(n-2)=\Phi(s,n)\geq\min\left\{\Phi(2,n),\Phi\left(\frac{n-1}{2},n\right)\right\}$ because the leading
coefficient of $\Phi(s,n)$ (view as a cubic polynomial of $s$) is positive, and $2\leq s\leq\frac{n-1}{2}$. According to $\frac{2}{3}<\alpha<1$
and $n\geq f(\alpha)=\frac{4}{1-\alpha}$, we obtain
\begin{align*}
\Phi(2,n)=&(1-\alpha)n^{2}+(6\alpha^{2}-4\alpha-3)n-28\alpha^{2}+42\alpha-10\\
\geq&(1-\alpha)\left(\frac{4}{1-\alpha}\right)^{2}+(6\alpha^{2}-4\alpha-3)\left(\frac{4}{1-\alpha}\right)-28\alpha^{2}+42\alpha-10\\
=&\frac{2}{1-\alpha}(14\alpha^{3}-23\alpha^{2}+18\alpha-3)\\
>&0,
\end{align*}
and
\begin{align*}
\Phi\left(\frac{n-1}{2},n\right)=&\frac{1}{8}((\alpha^{2}-3\alpha+2)n^{3}-(5\alpha^{2}-19\alpha+16)n^{2}\\
&+(3\alpha^{2}-25\alpha+34)n+\alpha^{2}+\alpha-20)\\
\geq&\frac{1}{8}((\alpha^{2}-3\alpha+2)\left(\frac{4}{1-\alpha}\right)^{3}-(5\alpha^{2}-19\alpha+16)\left(\frac{4}{1-\alpha}\right)^{2}\\
&+(3\alpha^{2}-25\alpha+34)\left(\frac{4}{1-\alpha}\right)+\alpha^{2}+\alpha-20)\\
=&\frac{1}{8(1-\alpha)^{2}}(\alpha^{4}-13\alpha^{3}+11\alpha^{2}+45\alpha-12)\\
>&0.
\end{align*}
Consequently, we conclude $\varphi_{B_1}(n-2)\geq\min\left\{\Phi(2,n),\Phi\left(\frac{n-1}{2},n\right)\right\}>0$ for $2\leq s\leq\frac{n-1}{2}$.
As $\eta_2<n-2<\rho_{\alpha}(K_1\vee(K_{n-2}\cup K_1))$ (see \eqref{eq:3.5}), we have
$$
\rho_{\alpha}(G_{s}^{2})<n-2<\rho_{\alpha}(K_1\vee(K_{n-2}\cup K_1))
$$
for $2\leq s\leq\frac{n-1}{2}$. This verifies Claim 1. \hfill $\Box$

In terms of \eqref{eq:3.1}, \eqref{eq:3.2} and Claim 1, we obtain
$$
\rho_{\alpha}(G)\leq\rho_{\alpha}(K_1\vee(K_{n-2}\cup K_1))
$$
with equality if and only if $G=K_1\vee(K_{n-2}\cup K_1)$. Which is a contradiction to the $A_{\alpha}$-spectral condition of Theorem 1.1. This
completes the proof of Theorem 1.1. \hfill $\Box$

\section{The proof of Theorem 1.2}

In this section, we prove Theorem 1.2, which provides an $A_{\alpha}$-spectral radius condition for a graph to be $t$-tough.

\medskip

\noindent{\it Proof of Theorem 1.2.} Suppose to the contrary that $G$ is not a $t$-tough graph, there exists some nonempty subset $S\subseteq V(G)$
such that $tc(G-S)-1\geq|S|$. Let $|S|=s$ and $c(G-S)=c$, then $tc-1\geq s$. If $n\geq(t+1)c-1$, then $G$ is a spanning subgraph of $G_1=K_{tc-1}\vee(K_{n_1}\cup K_{n_2}\cup\cdots\cup K_{n_c})$
for some positive integers $n_1\geq n_2\geq\cdots\geq n_c$ with $\sum\limits_{i=1}^{c}n_i=n-tc+1$. In terms of Lemma 2.2, we conclude
\begin{align}\label{eq:4.1}
\rho_{\alpha}(G)\leq\rho_{\alpha}(G_1)
\end{align}
with equality if and only if $G=G_1$.

Let $G_2=K_{tc-1}\vee(K_{n-(t+1)c+2}\cup(c-1)K_1)$. By virtue of Lemma 2.3, we obtain
\begin{align}\label{eq:4.2}
\rho_{\alpha}(G_1)\leq\rho_{\alpha}(G_2),
\end{align}
where the equality holds if and only if $(n_1,n_2,\ldots,n_c)=(n-(t+1)c+2,1,\ldots,1)$.

If $c=2$, then we have $G_2=K_{2t-1}\vee(K_{n-2t}\cup K_1)$. According to \eqref{eq:4.1} and \eqref{eq:4.2}, we get
$$
\rho_{\alpha}(G)\leq\rho_{\alpha}(K_{2t-1}\vee(K_{n-2t}\cup K_1)),
$$
where the equality holds if and only if $G=K_{2t-1}\vee(K_{n-2t}\cup K_1)$. In what follows, we consider $c\geq3$.

Recall that $G_2=K_{tc-1}\vee(K_{n-(t+1)c+2}\cup(c-1)K_1)$. According to $\alpha\in[\frac{1}{2},\frac{3}{4})$ and Lemma 2.5, we obtain
\begin{align}\label{eq:4.3}
\rho_{\alpha}(G_2)\leq&\frac{2e(G_2)(1-\alpha)}{n-1}+\alpha n-1\nonumber\\
=&\frac{(1-\alpha)((n-c+1)(n-c)+2(tc-1)(c-1))}{n-1}+\alpha n-1\nonumber\\
=&\frac{1-\alpha}{n-1}((2t+1)c^{2}-(2n+2t+3)c+n^{2}+n+2)+\alpha n-1.
\end{align}
Let $\varphi(c)=(2t+1)c^{2}-(2n+2t+3)c+n^{2}+n+2$. Note that $n\geq(t+1)c-1$. We conclude $3\leq c\leq\frac{n+1}{t+1}$. By a direct computation,
it follows from $n\geq\max\{5t^{2}+10t+1,\frac{12t(1-\alpha)-2\alpha+1}{3-4\alpha}\}\geq5t^{2}+10t+1>4t^{2}+6t+1$ that
$$
\varphi(3)-\varphi\left(\frac{n+1}{t+1}\right)=\frac{(n-3t-2)(n-4t^{2}-6t-1)}{(t+1)^{2}}>0,
$$
which implies that $\varphi(c)$ attains its maximum value at $c=3$ when $3\leq c\leq\frac{n+1}{t+1}$. Combining this with \eqref{eq:4.3},
$\alpha\in[\frac{1}{2},\frac{3}{4})$ and $n\geq\max\{5t^{2}+10t+1,\frac{12t(1-\alpha)-2\alpha+1}{3-4\alpha}\}\geq\frac{12t(1-\alpha)-2\alpha+1}{3-4\alpha}$,
we get
\begin{align}\label{eq:4.4}
\rho_{\alpha}(G_2)\leq&\frac{(1-\alpha)\varphi(3)}{n-1}+\alpha n-1\nonumber\\
=&\frac{(1-\alpha)(n^{2}-5n+12t+2)}{n-1}+\alpha n-1\nonumber\\
=&n-2+\frac{-(3-4\alpha)n+12t(1-\alpha)-2\alpha+1}{n-1}\nonumber\\
\leq&n-2+\frac{-(3-4\alpha)\cdot\frac{12t(1-\alpha)-2\alpha+1}{3-4\alpha}+12t(1-\alpha)-2\alpha+1}{n-1}\nonumber\\
=&n-2.
\end{align}

Since $K_{n-1}$ is a proper subgraph of $K_{2t-1}\vee(K_{n-2t}\cup K_1)$, we conclude
\begin{align}\label{eq:4.5}
\rho_{\alpha}(K_{2t-1}\vee(K_{n-2t}\cup K_1))>\rho_{\alpha}(K_{n-1})=n-2
\end{align}
by Lemmas 2.1 and 2.2. It follows from \eqref{eq:4.1}, \eqref{eq:4.2}, \eqref{eq:4.4} and \eqref{eq:4.5} that
$$
\rho_{\alpha}(G)\leq\rho_{\alpha}(G_1)\leq\rho_{\alpha}(G_2)\leq n-2<\rho_{\alpha}(K_{2t-1}\vee(K_{n-2t}\cup K_1)).
$$

In conclusion, we have
$$
\rho_{\alpha}(G)\leq\rho_{\alpha}(K_{2t-1}\vee(K_{n-2t}\cup K_1))
$$
with equality if and only if $G=K_{2t-1}\vee(K_{n-2t}\cup K_1)$. Which is a contradiction to the $A_{\alpha}$-spectral radius condition of
Theorem 1.2.

Let $n\leq(t+1)c-2$, then $G$ is a spanning subgraph of $K_{n-c}\vee cK_1$ with $c\geq\lceil\frac{n+2}{t+1}\rceil$. Let $G_3=K_{n-\lceil\frac{n+2}{t+1}\rceil}\vee\lceil\frac{n+2}{t+1}\rceil K_1$. By virtue of
Lemma 2.2, we conclude
\begin{align}\label{eq:4.6}
\rho_{\alpha}(G)\leq\rho_{\alpha}(K_{n-c}\vee cK_1)\leq\rho_{\alpha}(G_3),
\end{align}
with equalities if and only if $G=G_3$. Note that $2e(G_3)=(n-\lceil\frac{n+2}{t+1}\rceil)(n-\lceil\frac{n+2}{t+1}\rceil-1)+2\lceil\frac{n+2}{t+1}\rceil (n-\lceil\frac{n+2}{t+1}\rceil)
=(n-\lceil\frac{n+2}{t+1}\rceil)(n+\lceil\frac{n+2}{t+1}\rceil-1)<(n-\frac{n+2}{t+1})(n+\frac{n+2}{t+1})$. In view of Lemma 2.5, we possess
\begin{align}\label{eq:4.7}
\rho_{\alpha}(G_3)\leq&\frac{2e(G_3)(1-\alpha)}{n-1}+\alpha n-1\nonumber\\
<&\frac{(1-\alpha)(n-\frac{n+2}{t+1})(n+\frac{n+2}{t+1})}{n-1}+\alpha n-1\nonumber\\
=&n-2+\frac{\psi(n)}{(t+1)^{2}(n-1)},
\end{align}
where $\psi(n)=-(1-\alpha)n^{2}+((2-\alpha)t^{2}+(4-2\alpha)t+3\alpha-2)n-t^{2}-2t-5+4\alpha$. Notice that
$$
\frac{(2-\alpha)t^{2}+(4-2\alpha)t+3\alpha-2}{2(1-\alpha)}<5t^{2}+10t+1\leq n
$$
by $\alpha\in[\frac{1}{2},\frac{3}{4})$, $t\geq1$ and $n\geq\max\{5t^{2}+10t+1,\frac{12t(1-\alpha)-2\alpha+1}{3-4\alpha}\}$. For $t\geq2$, we have
\begin{align*}
\psi(n)\leq&\psi(5t^{2}+10t+1)\\
=&-(1-\alpha)(5t^{2}+10t+1)^{2}+((2-\alpha)t^{2}+(4-2\alpha)t+3\alpha-2)(5t^{2}+10t+1)\\
&-t^{2}-2t-5+4\alpha\\
=&(5t^{2}+10t+1)((4\alpha-3)t^{2}+(8\alpha-6)t+4\alpha-3)-t^{2}-2t-5+4\alpha\\
\leq&(5t^{2}+10t+1)(4(4\alpha-3)+2(8\alpha-6)+4\alpha-3)\\
&-t^{2}-2t-5+4\alpha \ \ \ \ \ \ \ \ \ \ \left(\mbox{since} \ t\geq2 \ \mbox{and} \ \alpha<\frac{3}{4}\right)\\
=&9(4\alpha-3)(5t^{2}+10t+1)-t^{2}-2t-5+4\alpha\\
<&0 \ \ \ \ \ \ \ \ \ \ \left(\mbox{since} \ t\geq2 \ \mbox{and} \ \alpha<\frac{3}{4}\right).
\end{align*}
We check that $\psi(n)<0$ also holds for $t=1$. Combining these with \eqref{eq:4.5}, \eqref{eq:4.6} and \eqref{eq:4.7}, we conclude
$$
\rho_{\alpha}(G)\leq\rho_{\alpha}(G_3)<n-2<\rho_{\alpha}(K_{2t-1}\vee(K_{n-2t}\cup K_1)),
$$
which contradicts $\rho_{\alpha}(G)\geq\rho_{\alpha}(K_{2t-1}\vee(K_{n-2t}\cup K_1))$. This completes the proof of Theorem 1.2. \hfill $\Box$

\section*{Data availability statement}

My manuscript has no associated data.

\section*{Declaration of competing interest}

The authors declare that they have no conflicts of interest to this work.

\section*{Acknowledgments}

The authors would like to express their sincere gratitude to the anonymous referees for their constructive corrections and valuable comments on this paper, which improved the presentation of this paper.
This work was supported by the Natural Science Foundation of Jiangsu Province (Grant No. BK20241949). Project ZR2023MA078 supported by Shandong Provincial Natural Science Foundation.

\end{document}